\documentclass[12pt,tbtags,leqno]{amsart}
\usepackage{amssymb}

\usepackage{amsmath,amsthm,bbm}
\usepackage{version}
\usepackage{xcolor,soul}

\usepackage{ucs}
\usepackage{mathrsfs}

\definecolor{dark_purple}{rgb}{0.4, 0.0, 0.4}
\definecolor{dark_green}{rgb}{0.0, 0.7, 0.0}

\newcommand\be{\begin{equation}}
\newcommand\ee{\end{equation}}
\newcommand\bt{\begin{theorem}}
\newcommand\et{\end{theorem}}
\newcommand\bprop{\begin{proposition}}
\newcommand\eprop{\end{proposition}}
\newcommand\beq{\begin{eqnarray*}}
\newcommand\eeq{\end{eqnarray*}}
\newcommand\bp{{\sc Proof: }}
\newcommand\ep{{}{\hfill $\Box$} \vskip 5pt \par}

\newcommand\bs{\vskip 12pt}

\newcommand\mccarthy{M\raise.45ex\hbox{c}Carthy}

\newcommand\C{\mathbb C}

\newcommand\T{\mathbb T}
\newcommand\D{\mathbb D}

\newcommand\Z{\mathbb Z}
\newcommand\NN{\mathbb N}

\newcommand\HH{{\mathbb H}}
\newcommand\RHP{{\mathbb P}}

\renewcommand\i{\infty}

\renewcommand\={\ = \ }
  \renewcommand\:[1]{\ {:=} \ }

\newcommand\vare{\varepsilon}
\newcommand\varphin{\phi^{[n]}}

\newcommand\phinp{\phi^{[n+1]}}
\newcommand\U{{\mathcal U}}
\newcommand\ZZ{{\mathcal Z}}
\newcommand\QQ{{\mathcal Q}}

\def\Id{\mathop{\rm Id}\nolimits}

\numberwithin{equation}{section}

\title{A note on composition operators on model spaces}

\author[I. Chalendar]{Isabelle Chalendar}
\address{I. Chalendar: Universit\'e Gustave Eiffel, LAMA, (UMR 8050), UPEM, UPEC, CNRS, F-77454, Marne-la-Vall\'ee, France}
\email{isabelle.chalendar@univ-eiffel.fr}

\author[P. Gumenyuk]{Pavel Gumenyuk}
\address{P. Gumenyuk: Department of Mathematics, Politecnico di Milano, via E. Bonardi 9, 20133 Milan, Italy.}
\email{pavel.gumenyuk@polimi.it}

\author[J.E. M\raise.5ex\hbox{c}Carthy]{John E. M\raise.5ex\hbox{c}Carthy}
\address{J.E. M\raise.5ex\hbox{c}Carthy: Department of Mathematics, Washington University in St. Louis, One Brookings Drive,
	St. Louis, MO 63130, USA}
\email{mccarthy@wustl.edu}

\thanks{J.M. was partially supported by National Science Foundation Grant DMS 2054199.}

\date{\today}



\begin{document}

\bibliographystyle{plain}

\theoremstyle{definition}
\newtheorem{definition}[equation]{Definition}
\newtheorem{lemma}[equation]{Lemma}
\newtheorem{proposition}[equation]{Proposition}
\newtheorem{theorem}[equation]{Theorem}
\newtheorem{claim}[equation]{Claim}
\newtheorem{question}[equation]{Question}
\newtheorem{corollary}[equation]{Corollary}
\newtheorem{example}[equation]{Example}
\newtheorem{remark}[equation]{Remark}
\newtheorem{problem}[equation]{Problem}

\begin{abstract} Motivated by the study of composition operators on model spaces launched by Mashreghi and  Shabankha  we consider the following problem: for a given inner function $\phi\not\in\mathsf{Aut}(\D)$, find a non-constant inner function~$\Psi$ satisfying the functional equation $\Psi\circ\phi=\tau\Psi$, where $\tau$ is a unimodular constant. We prove that this problem has a solution if and only  if $\phi$ is of positive hyperbolic step. More precisely, if this condition holds, we show that there is an infinite Blaschke product~$B$ satisfying the equation  for~$\tau=1$.  If in addition, $\phi$ is parabolic, we prove that the problem has a solution~$\Psi$ for \textit{any} unimodular~$\tau$. Finally, we show that if $\phi$ is of zero hyperbolic step, then no non-constant Bloch function~$f$ and no unimodular constant~$\tau$ satisfy $f\circ\phi=\tau f$.
\end{abstract}

\maketitle

%

\section{Introduction}\label{sec:intro}


A popular and successful subject in operator theory is the study of composition operators on Banach spaces of analytic functions.  We refer to the monographs \cite{CM95} by Cowen and McCluer and \cite{Sh93} by Shapiro for a comprehensive
presentation. In these two books as well as in the vast literature on this subject, these operators are considered on Banach spaces of analytic functions in the unit disk~$\D:=\{z\colon |z|<1\}$ such as  the Hardy space $H^2(\D)$. As a
consequence of the Littlewood subordination principle, for all holomorphic self-maps $\phi$ of $\D$, the composition operator $C_\phi$  is linear and bounded on $H^2(\D)$.

The seminal Beurling's Theorem
describes all the closed invariant subspaces  of the forward shift $S$ on $H^2(\D)$. They have the form
$\Theta H^2(\D)$, where $\Theta$ is an inner function, that is a bounded analytic function whose radial limits  are of modulus one almost everywhere. The so-called ``model spaces" denoted by $K_\Theta$ are defined to be their orthogonal complement in $H^2(\D)$. In other words, for $\Theta$ an inner function,
$K_\Theta:=H^2(\D)\cap (\Theta H^2(\D))^\perp$.  The $ K_\Theta$'s are therefore the closed invariant subspaces of the backward shift $S^*$.   Such model spaces are   Hilbert spaces, whose dimension is infinite  if and only if $\Theta $ is not a finite Blaschke product. Operator Theory
on those subspaces has undergone a great development in the past twenty years with the study of the so-called Trucated Toeplitz Operators (TTO) \cite{S07, GMR16, FM16}.

Much less has been done on
composition operators on model spaces.
This study was initiated by Mashreghi and  Shabankha  in~2013~\hbox{\cite{MS14,MS13}.}

The complete characterization of compact composition operators on model spaces was obtained by
Lyubarskii and Malinnikova \cite[Thm.\,1]{LM} who  proved that the composition operator ${C_\phi:K_\Theta\to H^2(\D)}$ is compact if and only if
\begin{equation}\label{LM}
	\lim_{|z|\to1^-}N_\phi(z)\frac{1-|\Theta(z)|^2}{1-|z|^2}=0,
\end{equation}
where $N_\phi$ is the Nevanlinna counting function defined by
\[N_\phi(z)\,:=\!\!\!\sum_{a\in\D,\,\phi(a)=z}\!\!\!(1-|a|)\quad\hbox{\sl\small(counted with multiplicity)}\quad \text{when $z\in\phi(\D)$}\]
 and $N_\phi(z):=0$ if $z\not\in \phi(\D)$.

  It follows that composition operators induced by inner functions $\phi$ on an infinite dimensional model space $K_\Theta$ are never compact. Indeed, by \cite[p.\,187]{sha87}, when $\phi$ is inner,   $N_\phi(z)\approx 1-|z|^2$ as $|z|\to 1$. Moreover, since $\Theta$ is not a finite Blaschle product, by a result due to Fatou there exists $(z_n)_{n\in\NN}\subset \D$ such that $|z_n|\to 1$ and $|\Theta(z_n)|$ does not converge to $1$.  Therefore   $\lim_{n\to\infty}N_\phi(z_n)\frac{1-|\Theta(z_n)|^2}{1-|z_n|^2}\neq 0$.\smallskip

From now on, suppose that $\phi$ is not the identity map and that $\Theta$ is not a finite Blaschke product.  In   \cite[Thm.\,4.1]{MS14}
they proved that if $\phi$ and $\Theta$ are inner functions,   then $C_\phi K_\Theta \subset K_\Theta$  in exactly three distinct cases:

\begin{itemize}
\item[I)] $\phi$  is an elliptic automorphism  and then
\[\Theta(z)=z\tau_p(z)^m \Psi ((\tau_p(z))^n),\]
where $p$ is the fixed point of $\phi$, $\tau_p(z):=\frac{p-z}{1-\overline{p}z}$, $m\in \NN_0$, $n\geq 2$, and $\Psi$ is  inner and not a finite Blachke product;   \medskip

\item[II)] $\phi$ has its Denjoy--Wolff point  $\alpha\in \T\:=\partial\D$ and $\Theta=z\Psi(z)$ where $\Psi$ is an inner function, not a finite Blaschke  product, such that $\Psi(\phi(z))=\tau \Psi(z)$ for some constant $\tau$ of modulus one;\medskip

\item[III)]  $\phi$ has its Denjoy--Wolff point  $\alpha\in \T$ and
\[\Theta=\gamma z\Psi(z)\prod_{n\geq 0}w (\phi^{[n]} (z)),\]
where $\phi^{[n]}$ stands for the $n$-th iterate of~$\phi$, $w$ is an inner function such that the product is convergent,
and $\Psi$ is an inner function, not a finite Blaschke  product, such that $\Psi(\phi(z))=\tau \Psi(z)$ for some constant $\tau$ of modulus one.
\end{itemize}

In cases I and II, the smallest model space containing the range of $C_\phi$ is $K_\Theta$. In Case III the smallest model space containing the range of $C_\phi$ is $K_\mu$ (see \cite[Thm.\,2.1]{MS13}), with $\mu w=\Theta$,   which implies that $K_\mu$ is strictly included in $K_\Theta$.\smallskip

The aim of this note is to study the solutions of  the equation
\begin{equation}\label{eq:main}
	\Psi\circ \phi \= \tau \Psi
\end{equation}
whenever $\phi$ has its Denjoy--Wolff point on the unit circle and  $\Psi$ is not a finite Blaschke product.  This contribution shows that, surprisingly, the theory of composition operators on model spaces is richer that expected.


\smallskip

Recall that the degree of an inner function $\Psi$ is  defined to be infinite when $\Psi$ is not a finite Blaschke product and is $d\in\NN$ when $\Psi$ is a finite Blaschke product with $d$ zeroes (taking into account their multiplicity). It is not difficult to check that the degree of the composition of two inner functions is equal to the product of their degrees.  Therefore if
  $\Psi$ is a finite Blaschke product of finite degree $d\in \NN$, the existence of  an inner function $\phi$ and $\tau\in\T$ such that \eqref{eq:main} holds implies  that $\phi$ is an automorphism. \smallskip

Moreover, as already noticed in \cite{MS14}, it is not difficult to see that when  $\phi\neq \Id$ has a fixed point $\alpha\in\D$,
then
\[\Psi\circ \phi=\tau \Psi\]
for some non-constant inner function $\Psi$ and
$\tau\in\T$ implies that  $\phi$ is  an elliptic automorphism.
Indeed, comparing the first non-constant terms in the Taylor expansions of $\Psi\circ\phi$ and $\tau\Psi$ at~$\alpha$, we see that $|\phi'(\alpha)|=1$ and hence, by the Schwarz--Pick Lemma, $\phi$ must be an automorphism of~$\D$.

In other words, when $\phi$ is not an automorphism, the existence of a non-constant inner function $\Psi$ and $\tau\in\T$ such that \eqref{eq:main} holds implies that $\Psi$ has infinite degree and $\phi$ has no fixed point in $\D$.  We can now formulate the main question at the center of our investigation.

\begin{question}
	\label{q1}
	Suppose $\phi : \D \to \D$ is inner, not an automorphism and with no fixed point in $\D$. Does there exist an inner function $\Psi$  satisfying
	\be
	\label{eq0}
	\Psi \circ \phi \= \tau \Psi
	\ee
	for some unimodular constant $\tau$?
\end{question}
We show that  examples can be constructed if and only if
$\phi$ is not of zero hyperbolic step.  \medskip

The paper is organized as follows. In Section~\ref{sec:2} we detail the classification of self-maps of the unit disk,  recalling the notion of positive and zero hyperbolic step, as well as properties of Abel's function and Julia's Lemma.

In Section~\ref{sec:3}  we prove our main result, Theorem~\ref{thmph1},
 asserting that if $\phi$ is an inner function with positive hyperbolic step, then there exists an infinite Blaschke product $B$ and  $\tau\in\T$ such that $B\circ \phi=\tau B$.

   In Section~\ref{sec:4} we strenghten Theorem \ref{thmph1} to show that $\tau$ can always be taken equal to $1$. Moreover, if
 $\phi$ is a parabolic map
  with positive hyperbolic step, we show that
  for every unimodular $\tau$  there exists an inner function $\Psi$ satisfying \eqref{eq:main}.
  	
 Section~\ref{sec:5} is devoted to the case when $\phi$ is a zero hyperbolic step map. In this case  we can answer negatively  our main question, and more than that, in Theorem~\ref{thmpzs}, we prove that  there does not exist a non-constant Bloch function and a unimodular number $\tau$
 satisfying $ f \circ \phi \= \tau f $. We conclude the paper with explicit examples of $\phi$
 of degree 2, $\Psi$ a Blaschke product, and $\tau$ satisfying \eqref{eq0}.

\section{Background}\label{sec:2}
We let $\D(z_0,r)$ denote the open unit disk in $\C$ centered at~$z_0$ of radius ${r>0}$, $\D:=\D(0,1)$, and let
$\HH$ denote the upper half-plane. For
$a \in \D$, we let $m_a$ be the Mobius automorphism
\[
m_a(z) := -\frac{a}{|a|}\frac{z-a}{1-\bar a z}~\text{~when~$~a\neq0$,\hskip1em and\hskip1em} m_0:=\Id.
\]
Denote by $\rho(\cdot,\cdot)$ the pseudo-hyperbolic distance in~$\D$, that is
$$
 \rho(z_1,z_2) \:= \left|\frac{z_2-z_1}{1-\overline{z_2}\,z_1}\right|,\quad z_1,z_2\in\D.
$$
Let $\phi: \D \to \D$ be a holomorphic self-map. We say it is of {\em positive hyperbolic step}
if and only if for some $z_0 \in \D$, we have
\be
\label{eqb1}
\lim_{n \to \i} \rho \big(\varphin(z_0), \phinp(z_0)\big) > 0.
\ee
Note that the above limit always exists because by the Schwarz--Pick lemma, the sequence is non-increasing. Moreover, if \eqref{eqb1} holds for some $z_0 \in \D$, then it holds for every $z_0 \in \D$; see e.g. \cite[Cor.\,4.6.9]{ab23}.

If $\phi$ is not of positive hyperbolic step, we say it is of {\em zero hyperbolic step.}
All elliptic maps are of zero hyperbolic step, hyperbolic maps are of positive hyperbolic step,
and parabolic maps can be either \cite{po79}.\smallskip

Our main result is that the answer to Question~\ref{q1} is yes if and only if $\phi$ is an inner function of positive hyperbolic step.\smallskip

Let $\phi$ be a non-elliptic self-map of $\D$.
Then Cowen \cite{cow81} proved that there is an Abel function, i.e. a holomorphic function
$h : \D \to \C$ that is  a solution
of the equation
\be
\label{eqch1}
h \circ \phi \= h + 1
\ee
and that is univalent on some $\phi$-absorbing domain $\U$ in $\D$.
Moreover, we can arrange that
\be
\label{eqom}
\Omega \:= \bigcup_{n \in {\mathbb N}}\  h(\D) - n
\ee
is $\HH$ if $\phi$ is  parabolic of positive hyperbolic step,
 $\C$ if it is of zero hyperbolic step, and a horizontal strip $\{ a < {\rm Im} (z) < b \}$ if it is
 hyperbolic. The Abel function is then unique up to an additive constant.\smallskip

 For $z_0 \in \D$ we define the {\em grand orbit of $z_0$}, denoted $\ZZ_{z_0}$, by
 \be
 \notag
\ZZ_{z_0} \:=  \big\{\zeta\in\D: \exists\,n,m \in \NN_0~ \text{s.t.}~\phi^{[n]}(\zeta)=\phi^{[m]}(z_0) \big\}.
 \ee

 For $\omega \in \T, \ M > 0$, let
\beq
H(\omega, M) &\:=& \big\{ z : \frac{|z-\omega|^2}{1 - |z|^2} < M \big\}\\
&=& \D\,\left( \frac{1}{M+1} \omega,\, \frac{M}{M+1}\right).
\eeq
Such sets are usually called \textit{horodisks} at~$\omega$.

Julia's Lemma \cite{ju20} (or e.g. \cite[Thm.\,5.9]{amy20}) says that
if $\phi$ admits at~$\omega$ a n.-t.(\,=\,non-tangential) limit $\eta\in\T$ and finite angular derivative~$a$, then
$\phi(H(\omega, M)) \subseteq H(\eta, a M)$.

\section{Positive Hyperbolic Step}\label{sec:3}

Throughout this section, we let  $\phi$ be an inner function of positive hyperbolic step, with Abel function $h$ as in \eqref{eqch1}. We shall let $\QQ$ denote the exceptional set of $\phi$, that is
\[
\QQ \:= \{ a : m_a \circ \phi~\text{has a  singular factor}\} .
\]
By Frostman's theorem, the set $\QQ$ is of logarithmic capacity zero, and in particular it has zero
area measure. See \cite[Sec.\,II.6]{gar81} for more details.

\begin{lemma}
\label{lemph1}
For each $z_0 \in \D$, the grand orbit $\ZZ_{z_0}$ satisfies the Blaschke condition.
\end{lemma}
\bp
Let $h$ be the Abel map. We have
\[
h(\ZZ_{z_0} ) \ \subseteq\ \{ h(z_0) + n : n \in \Z \} .
\]
Since $\phi$ is of positive hyperbolic step, the function
\[  \exp(2 \pi ih(z)) - \exp(2\pi ih(0))
\]
is bounded and non-constant. As the zero set of this function contains $\ZZ_{z_0}$,
we conclude that the grand orbit of $z_0$ satisfies the Blaschke condition.
\ep

Two sets $\ZZ_{z_0}$ and $\ZZ_{z_1}$ intersect if and only if they coincide, and in that case $h(z_1) - h(z_0) \in \Z$.

\begin{lemma}
\label{lemph2}
There exists $z_0$ with $\ZZ_{z_0} \cap \big(\QQ \cup\{z\in\D:\phi'(z)=0\}\big)$ empty.
\end{lemma}
\bp
Let $\QQ'$ denote $\QQ\cup\{z\in\D:\phi'(z)=0\}$.  Note  that $h(\QQ')$ has zero area measure. Indeed, otherwise there exists some $r < 1$ so
that $h(\QQ' \cap \D(0,r))$ has positive area. But on $\overline{\D(0,r)}$ the function
 $h$ has bounded derivative, so maps sets of area $0$ to sets of area $0$.

If $\ZZ_{z} \cap \QQ'$ is non-empty then
\be
\label{eqhp2}
h(z) \in \bigcup_{k\in \Z} h(\QQ') + k .
\ee
As $\bigcup_{k\in \Z} h(\QQ') + k $ has measure $0$, it follows that
\[
\{ h(z) : \ZZ_{z} \cap \QQ' \neq \emptyset \}
\]
has measure zero. Since $h(\D)$ has non-zero measure, it follows that there exists
 $z_0$ with $\ZZ_{z_0} \cap \QQ'$ empty.
\ep

\begin{theorem}
\label{thmph1}
Let $\phi$ be an inner function with positive hyperbolic step.
Then there exists a Blaschke product $B$ and a unimodular constant $\tau$ so that
\be
\label{eq03}
B \circ \phi \= \tau B  .
\ee
\end{theorem}
\bp
Choose $z_0$ so that $\ZZ_{z_0} \cap \big(\QQ \cup\{z\in\D:\phi'(z)=0\}\big)$ is empty.
Let $B$ be the Blaschke product whose zero set is $\ZZ_{z_0}$, with all zeros simple, that is
$$
  B(z):=\prod_{a\in \ZZ_{z_0}} m_a(z).
$$
Then the zero-set of the inner function $B \circ \phi$ is also $\ZZ_{z_0}$, again with all zeros simple, so
\[
B \circ \phi \= \tau u B,
\]
where $\tau$ is unimodular, and $u$ is either $1$ or a singular inner function. We will show that there is no singular inner factor.

Note that
\[
B \circ \phi \= \prod_{a\in\ZZ_{z_0}} m_{a} \circ \phi .
\]
By the choice of $z_0$, each factor on the right-hand side is a Blaschke product times a scalar, so therefore
so is their product.
Hence $B \circ \phi$ has no singular factor.
\ep

\section{Strengthening of Theorem~\protect{\ref{thmph1}}}\label{sec:4}

We can strengthen Theorem \ref{thmph1} in two ways.

\begin{theorem}
\label{thmph2}
Let $\phi$ be a parabolic inner function with positive hyperbolic step.
Then for every unimodular $\tau$  there exists an inner function $u$   so that
\be
\label{eqpph1}
u \circ \phi \= \tau u  .
\ee
\end{theorem}

Moreover, in Theorem \ref{thmph1} we can always choose $\tau = 1$.
\begin{theorem}
\label{thmhy2}
Let $\phi$ be an inner function with positive hyperbolic step.
Then there exists a Blaschke product $B$ satisfying
\[
B \circ \phi \= B .
\]
\end{theorem}

Before proving these theorems, we need the
 following results. The first one follows from \cite[Prop.\,3.10]{arbr16}; we include a proof for this special case.
\begin{proposition}
\label{propch1}
Suppose $\phi$ is a non-elliptic self-map of $\D$ with Abel map $h$ \and let {$\lambda\in\C\setminus\{0\}$}.
Then any holomorphic function $F$ satisfying
\be
\label{eqch2}
F \circ \phi \= \lambda F
\ee
is of the form
$F = G \circ h$ for some holomorphic  function $G$ on $\Omega$ that satisfies
\be
\label{eqch3}
G(w + 1) = \lambda G(w).
\ee
\end{proposition}

\bp
Note first that $h(a) = h(b)$ if and only if $\phi^{[n]}(a) = \phi^{[n]} (b)$ for some ${n\in\NN}$.
Indeed, as
\[
h \circ \phi^{[n]} \= h + n ,
\]
if $\phi^{[n]}(a) = \phi^{[n]} (b)$ then
$h(a) = h(b)$.
Conversely, if $h(a) = h(b)$, let $n\in\NN$ be large enough that
$\varphin(a)$ and $\varphin(b)$ lie in $\U$. As
\[
h \circ \varphin(a) \= h \circ \varphin (b)
\]
and $h$ is univalent on $\U$, it follows that $\varphin(a) = \varphin (b)$.

Suppose that $F$ satisfies \eqref{eqch2}. Then $F$ is constant on level-sets of $h$, so there is a well-defined function $G = F \circ h^{-1}$ on $h(\D)$. Moreover, $G$ is holomorphic at
every point $w$ such that some pre-image of $w$ under $h$ is not a critical point of $h$.
In particular, $G$ is holomorphic in $h(\U)$ and satisfies there $G(w+1)=\lambda G(w)$. We extend $G$ holomorphically to all of
$$
 \Omega \= \bigcup_{n\in\NN}h(\D)-n \= \bigcup_{n\in\NN}h(\U)-n
$$
by setting $G(w):=\lambda^{-n}\, G(w+n)$ for all ${w\in h(\U)-n}$ and ${n\in\NN}$.
\ep

\begin{proposition}
\label{propch2}
Suppose $\phi$ is an inner function
with positive hyperbolic step, and  with Abel map $h: \D \to \Omega$, where $\Omega$ is as in \eqref{eqom}. At a.e. point of $\T$, the Abel map $h$ has
a non-tangential limit that is in~$\partial \Omega$.
\end{proposition}
\bp
As $h : \D \to \Omega$, by Fatou's theorem, a.e. in~$\T$ it has non-tangential limits
that lie in $\overline{\HH} \cup \{ \infty \}$ if $\Omega = \HH$, and in
$\Omega \cup \{ \pm \infty \}$ if $\Omega$ is a horizontal strip.
Moreover these limits must be finite a.e., which follows from either
the  F.~and~M.~Riesz  theorem \cite[Thm. 2.5]{cl66} or Privalov's theorem \cite[Thm. 8.1]{cl66}.

By Theorem \ref{thmph1}, we know that there is a Blaschke product $B$ so that
$B \circ \phi = \tau B$.
By Proposition \ref{propch1}, we have $B = G \circ h$, where $G$ is holomorphic in~$\Omega$ and  satisfies $G(w+n)=\tau^n G(w)$ for any ${w\in\Omega}$ and any~${n\in\NN}$.
Since ${G(h(\D))=B(\D)\subset \D}$, and since for any $w\in\Omega$ there exists ${n\in\NN}$ such that $w+n\in h(\D)$, it follows that ${G(\Omega)\subset\D}.$
Therefore, if at some $\zeta\in\T$, the n.t. limit $h(\zeta)$ of~$h$ exists and belongs to~$\Omega$, then the n.t. limit of~$B$ at~$\zeta$ is $G(h(\zeta))\in\D$. Since~$B$ is inner, the latter may occur only on a set $E\subset\T$ of linear measure zero.
\ep

{\sc Proof of Theorem \ref{thmph2}:}
By Proposition \ref{propch2}, we know that $h$ has real n.t. boundary limits a.e.
Define for $0 \leq \theta \leq 2 \pi$,
\[
u_\theta \= \exp ( i\theta  h ) .
\]
Then each function $u_\theta$ is a bounded non-vanishing function on $\D$.
Moreover, $u_\theta$ has a unimodular boundary value a.e., so is inner.
The Abel equation implies
\be
\label{eqpph2}
\notag
u_\theta( \phi (z)) \= e^{ i \theta} u_\theta (z) .
\ee
This shows that we can solve \eqref{eqpph1} for every unimodular $\tau$.
\ep

{\sc Proof of Theorem \ref{thmhy2}:}
First we claim there is always an inner function $u_0$ satisfying
\be
\label{eqhy3}
u_0 (\phi(z)) \= u_0 (z) .
\ee
Indeed, if $\phi$ is parabolic, this was proved in Theorem \ref{thmph2}.
Assume instead that $\phi$ is hyperbolic, and define
$v := \exp( i h)$. Then $v$ is a bounded non-vanishing function that maps
$\D$ into an annulus,
satisfies $v (\phi(z)) = v(z)$, and by Proposition \ref{propch2}
 has non-tangential boundary limits on the boundary of the annulus a.e.

Let $f$ be an Ahlfors map from this annulus to the unit disk, see e.g. \cite[Ch.\,13]{bell16}; this is a two-to-one map
that maps the boundary of the annulus to the boundary of the disk.
Then
\[
u_0 \:= f \circ \exp(ih)
\]
is an inner function satisfying \eqref{eqhy3}.

Having found  an inner function satisfying \eqref{eqhy3},  we use
a  common trick (see for example \cite{s00,cc00,cgp15}) to extract a Blaschke product
that also satisfies \eqref{eqhy3}.
For each $a \in \D$, consider the function $\Psi_a:={m_a\circ u_0}$.
By Frostman's theorem \cite[Thm.\,II.6.4]{gar81}, for all~$a$ except for an exceptional set of capacity zero, this transformed function will be a Blaschke product times a unimodular scalar, and it is immediate that each $\Psi_a$
is also invariant under composition with~$\phi$.
\ep

\section{Zero hyperbolic step}\label{sec:5}
\label{secpar}

If $\phi$ is of zero hyperbolic step, it is either elliptic or parabolic. As explained in Section~\ref{sec:intro}, in the elliptic case, when the Denjoy--Wolff point is in $\D$, there can be no
non-constant holomorphic solution $f$ to the equation
\[
f \circ \phi \= \tau f
\]
with unimodular $\tau$ unless $\phi$ is an automorphism.

A holomorphic map $\phi: \D \to \D$ is parabolic if its Denjoy--Wolff point is on the boundary, and it has angular derivative $1$ there. Pommerenke~\cite{po79} showed that there exist parabolic self-maps of zero hyperbolic step as well as parabolic self-maps of positive hyperbolic step, see the definition in Section~\ref{sec:2}.

For convenience, we will change variables to the right-half plane $\RHP$ and assume the Denjoy--Wolff point is at $\i$.
Let us start with some initial point $z_0 = 1$, and define ${z_n := x_n + iy_n = \phi^{[n]}(z_0)}$.
Pommerenke proved the following \cite{po79}.

\bt
\label{thmpom}
Let $\phi$ be a parabolic self-map of $\RHP$. Define
\[
g_n(z) \= \frac{\varphin(z) - iy_n}{x_n} .
\]
Then $\lim_{n \to \i} g_n(z) = g(z)$ exists locally uniformly in $\RHP$,
and satisfies $g(\phi(z)) = \psi(g(z))$, where $\psi$ is a Moebius transformation of $\RHP$ that leaves $\i$ fixed.
Moreover, if $\phi$ is of positive hyperbolic step, then $\psi$ is parabolic, and if $\phi$ is of zero hyperbolic step
then $g(z) \equiv 1$.
\et

In \cite{bp79}, Baker and Pommerenke proved the following.
\bt
\label{thmbp}
Let $\phi$ be a parabolic  self-map of $\RHP$ of zero hyperbolic step
with Denjoy--Wolff point at $\i$. Define
\[
h_n(z) \= \frac{\phi^{[n]}(z) - z_n}{z_{n+1} - z_n} .
\]
Then $h(z) = \lim_{n \to \i} h_n(z)$ exists locally uniformly in $\RHP$ and satisfies
$h(\phi(z)) = z+1$.
\et

Using these results, we can show that $C_\phi$ can never have a Blaschke product as an eigenvector
with unimodular eigenvalue if $\phi$ is of zero hyperbolic step.
Indeed, we show slightly more. The Bloch space is the set of holomorphic functions $f$ on the unit disk
for which
\[
\sup_{z \in \D} (1 - |z|^2) | f'(z)| \ < \ \i.
\]
Note that it  follows from the Schwarz-Pick lemma that any bounded holomorphic function in~$\D$ is a Bloch function.

\bt
\label{thmpzs}
Let $\phi$ be a parabolic  self-map of~$\D$ with zero hyperbolic step.
Then there does not exist a non-constant Bloch function
 $f$
 and a unimodular number $\tau$
satisfying
\[
f \circ \phi \= \tau f .
\]
\et
We will use the following special case of \cite[Lemma~3.16]{arbr16}; an  elementary proof is included below.
\begin{lemma}
\label{leme1}
Let $\phi$ be a parabolic
self-map of $\D$ of zero hyperbolic step. For any pair of points $z_0, w_0 \in \D$,
we have
\[
\lim_{n \to \i} \rho\big(\varphin(z_0), \varphin(w_0)\big) \= 0.
\]
\end{lemma}
\bp
Let us change to the right-half plane, and assume that the Denjoy--Wolff point is at $\i$.
We can assume that $z_0 = 1$.
Write $z_n = x_n + iy_n= \varphin(z_0), w_n = \varphin(w_0)$.
We wish to show that
\[
\lim_{n \to \i} \rho(w_n, z_n) \= \lim_{n \to \i}\left|  \frac{w_n - z_n}{w_n + \overline{z_n}} \right| \= 0.
\]
Write
\be
\label{eqe1}
\frac{w_n - z_n}{w_n + \overline{z_n}}
\= \frac{w_n - z_n}{z_{n+1} - z_n} \frac{z_{n+1} - z_n}{w_n + \overline{z_n}} .
\ee
By Theorem \ref{thmbp}, the first fraction on the right-hand side of \eqref{eqe1} tends to $h(w_0)$.
We wish to show that the second fraction tends to $0$.
Let $0 < \vare < 1$. By two applications of Theorem \ref{thmpom} we get that there exists $N$ so that for $n \geq N$ we have
\beq
|w_n - z_n| &\ \leq \ & \vare x_n \\
|z_{n+1} - z_n| & \leq & \vare x_n .
\eeq
Therefore
\[
\left| \frac{z_{n+1} - z_n}{w_n + \overline{z_n}} \right|
\ ~\leq~ \frac{ \vare x_n}{(2 -  \vare ) x_n } .
\]
As $\vare$ is arbitrary, the limit of \eqref{eqe1} is $0$.
\ep

{\sc Proof of Theorem \ref{thmpzs}:}
Suppose $f$ is a Bloch function satisfying $ f \circ \phi = \tau f$, and let
\[
M \:= \sup_{z \in \D} (1 - |z|^2) | f'(z)| .
\]
Then for any points $z,w$ in the disk, we have the inequality
\[
| f(z) - f(w) | \ \leq \ M d_h(z,w),
\]
where
\[
d_h(z,w) \:= \log \frac{1 + \rho(z,w)}{1 - \rho(z,w)}
\]
is the hyperbolic metric. As $\tau $ is unimodular, we get
\beq
\left| \frac{f(z) - f(w)}{z-w} \right|
&\ \leq \ &
\left| \frac{f(\varphin(z)) - f(\varphin(w))}{z-w} \right| \\
&\leq& \frac{M}{|z-w|}\, d_h \big(\varphin(z), \varphin(w)\big) .
\eeq
As $d_h\big(\varphin(z), \varphin(w)\big) \to 0$ by Lemma \ref{leme1},
we conclude that $f$ is constant.
\ep

\section{Examples}\label{sec:6}

\begin{example}
Let \beq
b(z) &\:=& \frac{z + \alpha}{1 + \alpha z} \\
\phi(z) &\:=& b(z)^2 \=  \left( \frac{z + \alpha}{1 + \alpha z} \right)^2 ,
\eeq
where $ \frac 13 < \alpha < 1 $.
Then there exists a Blaschke product $B$ so that $B \circ \phi = - B$,
and hence $B^2 \circ \phi = B^2$.

\bp
We use the following properties of $\phi$:
\begin{enumerate}
\item
$\phi$ has Denjoy--Wolff point at $1 \in \partial \D$.
\item
The angular derivative at the Denjoy--Wolff point is
\[
a \= 2\, \frac{1-\alpha}{1+\alpha} \ < 1 .
\]
\item
$\phi$ is a finite Blaschke product, so its exceptional set is empty.
\item
The only zero of $\phi^\prime$ is at $-\alpha$.
\end{enumerate}
Let us divide the grand orbit $\ZZ_0$ into two sets:
\begin{eqnarray*}
{\mathcal M}_2 &\ :=\ & \{ z \in \ZZ_0: \ \exists\ n \geq 0 {\rm \ s.t.\ } \phi^{[n]}(z) = - \alpha \}
\\
{\mathcal M}_1 &\:= \ & \ZZ_0 \setminus {\mathcal M}_2 .
\end{eqnarray*}
We shall let $B$ be the Blaschke product with zeroes in $\ZZ_0$,
with multiplicity $2$ for points in ${\mathcal M}_2$ and multiplicity $1$ for points in ${\mathcal M}_1$.

The point ${z_0:=0}$ belongs to the boundary of the horodisk $H(1,1)$. Hence by Julia's Lemma, see Section~\ref{sec:2}, the pre-images of~$z_0$ lie outside $H(1,\frac 1a)$.
Their pre-images lie outside $H(1,a^{-2})$, and so on.
For $m\in\NN$, let ${z_m}$ denote $\phi^{[m]}(z_0)$.
Consider the pre-images of $z_1$. There are two: $z_0$,  and the solution to
$b(z) = - b(z_0)$.
This point is \[
\zeta_0 \ := \ - \frac{2\alpha}{1+\alpha^2},
\]
 which is in the boundary of $H(1,  \frac{1-\zeta_0}{1+ \zeta_0} ) = H(1,\frac{4}{a^2})$.

For each $m\in\NN$, consider the pre-images of $\phi(z_m)$.
There are two. One is $z_{m-1}$; call the other $\zeta_{m-1}$. This is the solution to the equation
\be
\label{eq4}
b(\zeta_{m-1}) \= - b(z_{m-1}).
\ee
This gives
\be
\label{eq5}
\zeta_{m-1} \= \frac{\zeta_0 - z_{m-1} }{1 - \zeta_0 z_{m-1}}.
\ee
Notice that each $\zeta_{m-1}$ is negative.

 By an argument similar to the proof of Theorem \ref{thmph1}, we have
$B \circ \phi = \tau B$ for some $\tau\in\T$.
Indeed, $B \circ \phi$ is a unimodular multiple of a Blaschke product,
it has zeroes only in the set $\ZZ_0$, and their multiplicity is $1$ on
the set  ${\mathcal M}_1$ and
$2$
 on  ${\mathcal M}_2$.

Note that $B(0) = 0$.
All the zeroes of $B$ are symmetric with respect to the real axis.
 This is because $\phi^{[n]}(\bar z) = \overline{\phi^{[n]}(z)}$, so
 if $z \in \ZZ_0$, we have $\phi^{[n]} (z) = z_m$ for some $m, n\in\NN_0$, and
 as $z_m$ is always real, this means that $\phi^{[n]} (\bar z) = z_m$ also, so $\bar z \in
 \ZZ_0$.
Therefore $B$ is real on $(-1,1)$. Moreover, as $\phi$ is real on $(-1,1)$, so  is $B \circ \phi$. Therefore $\tau$ is real.

There are no zeroes of $B$ on the line segment between $0$ and $z_1 = \phi(0) = \alpha^2$.
Indeed, all the pre-images of $z_0$ are outside $H(1,1)$, so in particular do not lie in the set $(0,1)$.
Moreover,  each $\zeta_n$ is negative, so it and all its pre-images also lie outside $H(1,1)$.
So the only zeroes of $B$ that lie on the line segment $(0,1)$ are the points $z_m$, for $m \geq 1$. These points increase. Therefore there are no intermediate zeroes between $z_0$ and $z_1$.

As each zero of $B$ on~$[0,1)$ is of multiplicity $1$, the sign of $B'$ will alternate as one moves along the real axis.
By the chain rule,
\[
B'(z_1) \phi'(0) \= \tau B'(0) .
\]
As $\phi'(0) = 2 \alpha (1 - \alpha^2) > 0$,  and $B'(0)$ and $B'(z_1)$ are opposite signs, we conclude that $\tau$ must be negative.
\ep
\end{example}

\begin{example}
Let $\alpha = \frac 13$ in the previous example, so
\[
\phi(z) \ =\ \frac{(z+1/3)^2}{(1 + z/3)^2}.
\]
Then $\phi$ has a fixed point at $1$, and $\phi'(1) = 1$, so it is parabolic. We show that it is of zero hyperbolic step,
 and hence according to Theorem~\ref{thmpzs}, the answer to Question~\ref{q1} is negative for this inner function~$\phi$.

\bp
We calculate, if $z \in (-1,1)$ is real:
\beq
\rho(z,\phi(z) ) &\ = \ & \left| \frac{ z - \left( \frac{z+ \frac 13}{1 + \frac 13 z}\right)^2}{1 - z  \left( \frac{z+ \frac 13}{1 + \frac 13 z}\right)^2} \right| \\
&=& \frac{ (1-z)^2} {9z^2 + 14 z + 9} .
\eeq
If we let $z_0 = 0$, then $z_n \to 1$, and each $z_n$ is real.
But we see that as $z_n \to 1$, we have
$\rho(z_n, z_{n+1}) \to 0$.
\ep
\end{example}

\begin{remark}
Another (less direct) way to show that $\phi$ in the above example is of zero hyperbolic step is based on the fact that all orbits of a parabolic self-map of \textit{positive} hyperbolic step converge to the Denjoy\,--\,Wolff point \textit{tangentially} to~$\T$; see \cite[Rmk.\,1]{po79} or~\cite[Cor.\,4.6.10]{ab23}.
\end{remark}

\bs
\bs
\bs
{\bf Acknowledgement:} This work was done while the third author was visiting the Mathematics Department at Universit\'e Gustave Eiffel. He would like to thank the Department for its hospitality.

\bibliography{references3}

\begin{thebibliography}{10}

\bibitem{ab23}
Marco Abate.
\newblock {\em Holomorphic dynamics on hyperbolic {R}iemann surfaces},
  volume~89 of {\em De Gruyter Studies in Mathematics}.
\newblock De Gruyter, Berlin, 2023.

\bibitem{amy20}
Jim Agler, John~E. M\raise.45ex\hbox{c}Carthy, and N.~J. Young.
\newblock {\em Operator Analysis: Hilbert Space Methods in Complex Analysis}.
\newblock Cambridge Tracts in Mathematics. Cambridge University Press, 2020.

\bibitem{arbr16}
Leandro Arosio and Filippo Bracci.
\newblock Canonical models for holomorphic iteration.
\newblock {\em Trans. Amer. Math. Soc.}, 368(5):3305--3339, 2016.

\bibitem{bp79}
I.~N. Baker and Ch. Pommerenke.
\newblock On the iteration of analytic functions in a halfplane. {II}.
\newblock {\em J. London Math. Soc. (2)}, 20(2):255--258, 1979.

\bibitem{bell16}
Steven~R. Bell.
\newblock {\em The {C}auchy transform, potential theory and conformal mapping}.
\newblock Chapman \& Hall/CRC, Boca Raton, FL, second edition, 2016.

\bibitem{cc00}
G.~Cassier and I.~Chalendar.
\newblock The group of the invariants of a finite {Blaschke} product.
\newblock {\em Complex Variables, Theory Appl.}, 42(3):193--206, 2000.

\bibitem{cl66}
E.~F. Collingwood and A.~J. Lohwater.
\newblock {\em The theory of cluster sets}.
\newblock Cambridge Tracts in Mathematics and Mathematical Physics, No. 56.
  Cambridge University Press, Cambridge, 1966.

\bibitem{cow81}
Carl~C. Cowen.
\newblock Iteration and the solution of functional equations for functions
  analytic in the unit disk.
\newblock {\em Trans. Amer. Math. Soc.}, 265(1):69--95, 1981.

\bibitem{CM95}
Carl~C. Cowen and Barbara~D. MacCluer.
\newblock {\em Composition operators on spaces of analytic functions}.
\newblock Boca Raton, FL: CRC Press, 1995.

\bibitem{FM16}
Emmanuel Fricain and Javad Mashreghi.
\newblock {\em The theory of {{\(\mathcal{H}(b)\)}} spaces. {Volume} 1},
  volume~20 of {\em New Math. Monogr.}
\newblock Cambridge: Cambridge University Press, 2016.

\bibitem{GMR16}
Stephan~Ramon Garcia, Javad Mashreghi, and William~T. Ross.
\newblock {\em Introduction to model spaces and their operators}, volume 148 of
  {\em Camb. Stud. Adv. Math.}
\newblock Cambridge: Cambridge University Press, 2016.

\bibitem{gar81}
John~B. Garnett.
\newblock {\em Bounded Analytic Functions}.
\newblock Academic Press, New York, 1981.

\bibitem{cgp15}
P.~Gorkin I.~Chalendar and J.R. Partington.
\newblock Inner functions and operator theory.
\newblock {\em North-W. Eur. J. of Math.}, 1:9--28, 2015.

\bibitem{ju20}
Gaston Julia.
\newblock Extension nouvelle d'un lemme de {Schwarz}.
\newblock {\em Acta Math.}, 42:349--355, 1920.

\bibitem{LM}
Yurii~I. Lyubarskii and Eugenia Malinnikova.
\newblock Composition operators on model spaces.
\newblock In {\em Recent trends in analysis. Proceedings of the conference in
  honor of Nikolai Nikolski on the occasion of his 70th birthday, Bordeaux,
  France, August 31 -- September 2, 2011}, pages 149--157. Bucharest: The Theta
  Foundation, 2013.

\bibitem{MS14}
J.~Mashreghi and M.~Shabankhah.
\newblock Composition of inner functions.
\newblock {\em Can. J. Math.}, 66(2):387--399, 2014.

\bibitem{MS13}
Javad Mashreghi and Mahmood Shabankhah.
\newblock Composition operators on finite rank model subspaces.
\newblock {\em Glasg. Math. J.}, 55(1):69--83, 2013.

\bibitem{po79}
Ch. Pommerenke.
\newblock On the iteration of analytic functions in a halfplane.
\newblock {\em J. London Math. Soc. (2)}, 19(3):439--447, 1979.

\bibitem{S07}
Donald Sarason.
\newblock Algebraic properties of truncated {Toeplitz} operators.
\newblock {\em Oper. Matrices}, 1(4):491--526, 2007.

\bibitem{sha87}
Joel~H. Shapiro.
\newblock The essential norm of a composition operator.
\newblock {\em Ann. of Math. (2)}, 125(2):375--404, 1987.

\bibitem{Sh93}
Joel~H. Shapiro.
\newblock {\em Composition operators and classical function theory}.
\newblock Universitext. New York: Springer-Verlag, 1993.

\bibitem{s00}
Joel~H. Shapiro.
\newblock What do composition operators know about inner functions?
\newblock {\em Monatsh. Math.}, 130(1):57--70, 2000.

\end{thebibliography}
\end{document}